\documentclass[12pt,leqno]{amsart}
\usepackage{amsmath,amssymb,amsfonts}
\usepackage{eucal}

\newcommand \comment[1]{}			%  Silent version.
\newcommand \mylabel[1]{\label{#1}\comment{{\rm \{#1\} }}}

\newtheorem*{thm}{Theorem}
 
\renewcommand{\phi}{\varphi}

\newcommand\bbm{\begin{bmatrix}}
\newcommand\ebm{\end{bmatrix}}
 
\newcommand\cP{\mathcal{P}}
\newcommand\cS{\mathcal{S}}
\newcommand\bbN{\mathbb{N}}
\newcommand\bbZ{\mathbb{Z}}
\newcommand\fG{\mathfrak G}
\newcommand\fH{\mathfrak H}
\newcommand\fS{\mathfrak S}

\newcommand\D{\Delta}
\newcommand\G{\Gamma}
\newcommand\X{\mathbf X}
\newcommand\chiset{\chi^{\mathrm{set}}}

\newcommand\Fr{\mathrm{Fr}}
\newcommand \inv{^{-1}}
\newcommand\GF{\operatorname{GF}}

\newcommand\gauss[2]{\bbm #1\\ #2 \ebm}

\allowdisplaybreaks

\begin{document}

 \title{A new distribution problem of balls into urns, and how to color a graph by different-sized sets}
 \author{Thomas Zaslavsky}
 \address{Department of Mathematical Sciences\\ Binghamton University (SUNY)\\ Binghamton, NY 13902-6000 \\ U.S.A.}
\email{zaslav@math.binghamton.edu}

\begin{abstract}
Set-coloring a graph means giving each vertex a subset of a fixed color set so that no two adjacent subsets have the same cardinality.  When the graph is complete one gets a new distribution problem with an interesting generating function.  We explore examples and generalizations.
\end{abstract}

\date{5 July 2006; first version 25 June 2006.  This version \today}
 
 \maketitle

\subsection*{Balls into urns} We have $n$ labelled urns and an unlimited supply of balls of $k$ different colors.  Into each urn we want to put balls, no two the same color, so that the number of colors in every urn is different.  Balls of the same color are indistinguishable and we don't care if several are in an urn.  How many ways are there to do this?  
 (The reader will note the classical terminology.  Our question appears to be new but it could as easily have been posed a hundred years ago.)
 Call the answer $\chiset_n(k)$.  We  form the exponential generating function,
 \[
 \X(t) := \sum_{n=0}^\infty \chiset_n(k) \frac{t^n}{n!} \ ,
 \]
taking $\chiset_0(k) = 1$ in accordance with generally accepted counting principles.  Then we have the generating function formula
\begin{equation}\mylabel{E:urns}
\X(t) = \prod_{j=0}^k \Big[ 1 + \binom{k}{j} t \Big] .
\end{equation}
For the easy proof, think about how we would choose the sets of colors for the urns.  We pick a subset of $n$ integers, $\{j_1<j_2<\cdots j_n\} \subseteq \{0,1,\ldots,k\}$, and assign each integer to a different urn; then we choose a $j_i$-element subset of $[k] := \{1,2,\ldots,k\}$ for the $i$th urn.  The number of ways to do this is
\[
\sum_{S\subseteq[k]: |S|=n} n! \, \prod_{j\in S} \binom{k}{j}.
\]
Forming the exponential generating function, the rest is obvious.
 
There are several interesting features to the question and its answer.  
First of all, as far as I know the question is a new distribution problem.  
Second, the sequence $\chiset_0(k), \chiset_1(k), \ldots, \chiset_{k+1}(k)$, besides (obviously) being increasing, is logarithmically concave, because the zeros of its generating function are all negative real numbers.
Third, the theorem generalizes to graphs and can be proved by means of M\"obius inversion over the lattice of connected partitions of the vertex set, just as one proves the M\"obius-function formula for the chromatic polynomial.  
%Because of the graphical generalization we call it \emph{set coloring} (the graph in the balls-into-urns problem being the complete graph; see later).  
Fourth, \eqref{E:urns} and the graphical extension generalize to formulas in which the binomial coefficients are replaced by arbitrary quantities.  
Finally, this way of putting balls into urns, and its graphical generalization, are really a problem of coloring gain graphs by sets, which suggests a new kind of gain-graph coloring; we discuss this briefly at the end.
 
Some elementary notation: $\bbN$ denotes the set of nonnegative integers and $[n] := \{1,2,\ldots,n\}$ for $n\geq0$; $[0]$ is the empty set.  Furthermore, $\cP_k$ denotes the power set of $[k]$.

To set the stage for the graph theory, first we generalize Equation \eqref{E:urns}.  Let $\alpha := (\alpha_j)_0^\infty$ be a sequence of numbers, polynomials, power series, or any quantities for which the following expressions and generating functions, including those in Equation \eqref{E:gf}, are defined.  Let $\beta_r := \sum_{j=0}^\infty \alpha_j^r$.  Let $\chi_n(\alpha)$ := the sum of\/ $\prod_1^n \alpha_{f(i)}$ over all injective functions $f : [n] \to \bbN$.  Then, generalizing \eqref{E:urns}, and with a similar proof, we have
\begin{equation}\mylabel{E:gf}
\X_a(t) := \sum_{n=0} \chi_n(\alpha) \frac{t^n}{n!} =  \prod_{j=0}^\infty \big[ 1 + \alpha_j t \big] .
\end{equation}
 As with the set-coloring numbers, if $\alpha$ is nonnegative and is a finite sequence $(\alpha_j)_{j=0}^k$, then the sequence $(\chi_n(\alpha))$ is logarithmically concave.  We can even closely approximate the index $m$ of the largest $\chi_n(\alpha)$.  Darroch's Theorem 3 \cite{Darroch} says that $m$ is one of the two nearest integers to 
 \[
 M := k+1 - \sum_{j=0}^k \frac{1}{1+\alpha_j} ,
 \]
and $m=M$ if $M$ is an integer.\footnote{I thank Herb Wilf for telling me of Darroch's theorem and reminding me about logarithmic concavity.}

A combinatorial problem that falls under Equation \eqref{E:gf} is filling urns from the equivalence classes of a partition.  We have a finite set $\cS$ with a partition that has $k+1$ blocks $\cS_0, \cS_1,\ldots, \cS_k$.  
We want the number of ways to put one ball into each of $n$ labelled urns with no two from the same block.  Call this number $\chi_n(\pi)$.  The generating function is \eqref{E:gf} with $\alpha_j = |\cS_j|$.  
It is clear that $\chi_n(\pi)$ increases with its maximum at $n=k$.  
As an example let $\cS =$ the lattice of flats of a rank-$k$ matroid, two flats being equivalent if they have the same rank; then $\alpha_j = W_j$, the number of flats of rank $j$ (the Whitney number of the second kind).  In particular, if $\cS =$ the lattice of subspaces of the finite vector space $\GF(q)^k$, the rule being that each urn gets a vector space of a different dimension, then 
\[
\X_\cS(t) = \prod_{j=0}^k \left( 1 + \gauss{k}{j} t \right),
\]
where $\gauss{k}{j}$ is the Gaussian coefficient.  
For a similar example where the $\alpha_j$ are (the absolute values of) the Whitney numbers of the first kind, take $\cS$ to be the broken circuit complex of the matroid.

\subsection*{Graphs}  In the graphical generalization we have $\D$, a graph on vertex set $V=[n]$.  $\Pi(\D)$ is the set of \emph{connected partitions} of $\D$, that is, partitions $\pi$ of $V$ such that each block $B \in \pi$ induces a connected subgraph.  The set $\Pi(\D)$, ordered by refinement, is a geometric lattice with bottom element $\hat0$, the partition in which every block is a singleton.  
A \emph{set $k$-coloring} of $\D$ is a function $c: V \to \cP_k$ that assigns to each vertex a subset of $[k]$, and it is \emph{proper} if no two adjacent vertices have colors (that is, sets) of the same cardinality.  We define the \emph{set-coloring function} $\chiset_\D(k)$ to be the number of proper set $k$-colorings of $\D$.  This quantity is positive just when the chromatic number of $\D$ does not exceed $k+1$.

The \emph{extended Franel numbers} are
\[
\Fr(k,r) := \sum_{j=0}^k \binom{k}{j}^r 
\]
for $k, r \geq 0$.  (The Franel numbers themselves are the case $r=3$ \cite[Sequence A000172]{OEIS}.  There is a nice table of small values of the extended numbers at \cite{BinomMW}.  There are closed-form expressions when $r \leq 2$ but not otherwise.)  
The set-coloring function satisfies
\begin{equation}\mylabel{E:set-coloring}
\chiset_\D(k) = \sum_{\pi\in\Pi(\D)} \mu(\hat0,\pi) \prod_{B\in\pi} \Fr(k,|B|)
\end{equation}
where $\mu$ is the M\"obius function of $\Pi(\D)$.

It is amusing to see the high-powered machinery involved in deriving \eqref{E:urns} from \eqref{E:set-coloring}.  We outline the method.  
Obviously, $\chiset_{K_n}(k) = \chiset_n(k)$.  In \eqref{E:set-coloring} we substitute the known value $\mu(\hat0,\pi) = \prod_{B\in\pi} [ -(-1)^{|B|}(|B|-1)! ]$.  Then we apply the exponential formula to the exponential generating function, substituting $y = -\binom{k}{j}t$ in $\log(1-y) = -\sum_{n=1}^\infty y^n/n$ and finding that the exponential and the logarithm cancel.

Rather than proving Equation \eqref{E:set-coloring} itself, we generalize still further; the proof is no harder.  Define 
 \[
 \chi_\D(\alpha) := \sum_f  \prod_{i=1}^n \alpha_{f(i)},
 \]
summed over all functions $f: V \to \bbN$ such that $f(i) \neq f(j)$ if $i$ and $j$ are adjacent; that is, over all proper $\bbN$-colorings of $\D$.  One could think of $f$ as a proper $\bbN$-coloring weighted by $\prod \alpha_{f(i)}$.  (Again, we assume $\alpha$ has whatever properties are required to make the various sums and products in the theorem and its proof meaningful.  A sequence that is finitely nonzero will satisfy this requirement.) 
 
 \begin{thm}\mylabel{T:graph-labels}
We have
 \[
 \chi_\D(\alpha) = \sum_{\pi\in\Pi(\D)} \mu(\hat0,\pi) \prod_{B\in\pi} \beta_{|B|}.
 \]
\end{thm}
 
To derive Equation \eqref{E:set-coloring} we set $\alpha_j = \binom{k}{j}$.  It is easy to see that the left-hand side of the theorem equals $\chiset_\D(k)$.
 
 \begin{proof}
The method of proof is the standard one by M\"obius inversion.  For $\pi\in\Pi(\D)$ define 
\[
g(\pi) = \sum_f \prod_1^n \alpha_{f(i)},
\]
summed over functions $f: V \to \bbN$ that are constant on blocks of $\pi$, and 
\[
h(\pi) = \sum_f \prod_1^n \alpha_{f(i)},
\]
summed over every such function whose values differ on blocks $B, B' \in \pi$ that are joined by one or more edges.  It is clear that
\[
g(\pi') = \sum_{\pi\geq\pi'} h(\pi)
\]
for every $\pi'\in \Pi(\D)$, $\pi$ also ranging in $\Pi(\D)$.  By M\"obius inversion,
\begin{equation}\mylabel{E:mu}
h(\pi') = \sum_{\pi\geq\pi'} \mu(\pi',\pi) g(\pi).
\end{equation}
Set $\pi' = \hat0$ and observe that $h(\hat0) = \chi_\D(\alpha)$.

To complete the proof we need a direct calculation of $g(\pi)$.  We may choose $f_B \in \bbN$ for each block of $\pi$ and define $f(i)=f_B$ for every $i\in B$; then 
\[
g(\pi) = \prod_{B\in\pi} \sum_{j=0}^\infty \alpha_j^{|B|} = \prod_{B\in\pi} \beta_{|B|}.
\]
Combining with Equation \eqref{E:mu}, we have the theorem.
\end{proof}

As with our original balls-into-urns problem, there is a combinatorial special case where we color $\D$ from a set $\cS$ with a partition $\pi$, so that no two adjacent vertices have equivalent colors.  We call this \emph{coloring from a partitioned set} and denote the number of ways to do it by $\chi_\D(\pi)$.

\newcommand\colmin{.3in}
\begin{table}[htdp]
\begin{center}
\begin{tabular}{rc | rrr rrr rrr}
 & $k$ & 0 & 1 & 2 & 3 & 4 & 5 & 6 & 7  \\
$n$ && \hspace{\colmin} & \hspace{\colmin} & \hspace{\colmin} & \hspace{\colmin} & \hspace{\colmin} & \hspace{\colmin} & \hspace{\colmin} & \hspace{\colmin} & \hspace{\colmin} \\
\hline
0 && 1 & 1 & 1 & 1 & 1 		& 1 	& 1 	& 1 \\
1 && 1 & 2 & 4 & 8 & 16 	& 32 	& 64 	& 128 \\
2 && 0 & 2 & 10 & 44 & 186 	& 772 	& 3172 	& 12952 \\
3 && 0 & 0 & 12 & 144 & 1428 	& 13080 & 115104 & 989184 \\
4 && 0 & 0 & 0 & 216 & 6144 	& 139800 & 2821464 & 53500944 \\
5 && 0 & 0 & 0 & 0 & 11520 	& 780000 & 41472000 & 1870310400 \\
6 && 0 & 0 & 0 & 0 & 0 		& 1800000 & 293544000 & 37139820480 \\
7 && 0 & 0 & 0 & 0 & 0 		& 0 	& 816480000 & 325275955200 \\
8 && 0 & 0 & 0 & 0 & 0 		& 0 	& 0 	& 1067311728000 \\
9 && 0 & 0 & 0 & 0 & 0 		& 0 	& 0 	& 0 \\
%  &  &  &  &  &  &  &  &  &  \\       %SPARE BLANK LINE
\end{tabular}
\bigskip
\caption{Values of $\chiset_n(k)$ for small $n$ and $k$.}\mylabel{Tb:urns}
\end{center}

\end{table}

\subsection*{Examples}  The table %\myref{Tb:urns} 
shows some low values of $\chiset_n(k)$, and the list below has formulas for special cases.  
We also calculate two graphical set-chromatic functions. A trivial one is $\chiset_\D(k)$ for $\bar K_n$, the graph with no edges, since $\chiset_\D$ is multiplicative over connected components, and it is not hard (if tedious) to do graphs of order at most $3$, such as the $3$-vertex path $P_3$.  Here are some examples:
 \begin{align*}
\chiset_0(k) &= 1, \\
\chiset_1(k) &= 2^k, \\
\chiset_2(k) &= 2^{2k} - \binom{2k}{k}, \\
\chiset_3(k) &= 2^{3k} - 3\cdot2^k\binom{2k}{k} + 2\cdot\Fr(k,3) , \\
\chiset_n(k) &= 0 \text{ when } k < n-1, \\
\chiset_n(n-1) &= n! \, \binom n0 \binom n1 \cdots \binom nn , \\
\chiset_{P_3}(k) &= 2^{3k} - 2\cdot2^k\binom{2k}{k} + \Fr(k,3), \\
\chiset_{\bar K_n}(k) &= 2^{nk}.
\end{align*}
The table entries for $n>3$ were obtained from the preceding formulas and, with the help of Maple, from the generating function \eqref{E:urns}.

The table shows that the values of $\chiset_2(k)$ match the number of rooted, $k$-edge plane maps with two faces \cite[Sequence A068551]{OEIS}.  The two sequences have the same formula.  It would be interesting to find a bijection.  A casual search of \cite{OEIS} did not reveal any other known sequences in the table that were not obvious.

\subsection*{Gain graphs} Set coloring began with an idea about gain graph coloring when the gains are permutations of a finite set.

Take a graph $\G$, which may have loops and parallel edges, and assign to each oriented edge $e_{ij}$ an element $\phi(e_{ij})$ of the symmetric group $\fS_k$ acting on $[k]$, in such a way that reorienting the edge to the reverse direction inverts the group element; symbolically, $\phi(e_{ji}) = \phi(e_{ij})\inv$.  We call $\phi$ a \emph{gain function} on $\G$, and $(\G,\phi)$ is a \emph{permutation gain graph} with $\fS_k$ as its \emph{gain group}.  
A \emph{proper set coloring} of $(\G,\phi)$ is an assignment of a subset $S_i \subseteq [k]$ to each vertex $i$ so that for every oriented edge $e_{ij}$, $S_j \neq S_i \phi(e_{ij})$.  One way to form a permutation gain graph is to begin with a simple graph $\D$ on vertex set $[n]$ and replace each edge ${ij}$ by $k!$ edges $(g,{ij})$, each labelled by a different element $g$ of the gain group $\fS_k$.  (Then the notations $(g,{ij})$ and $(g\inv,{ji})$ denote the same edge.)  
We call this the \emph{$\fS_k$-expansion} of $\D$ and write it $\fS_k\D$.  Now a proper set coloring of $\fS_k\D$ is precisely a proper set coloring of $\D$ as we first defined it: an assignment to each vertex of a subset of $[k]$ so that no two adjacent vertices have sets of the same size.  Thus I came to think of set-coloring a graph.

Our calculations show that the number of proper set colorings of a graph $\D$, or equivalently of its $\fS_k$-expansion, is exponential in $k$.  There is a standard notion of coloring of a gain graph with gain group $\fG$, in which the colors belong to a group $\fH=\fG\times\bbZ_k$ and there is a chromatic function, a polynomial in $|\fH|$, that generalizes the chromatic polynomial of an ordinary graph and has many of the same properties, in particular satisfying the deletion-contraction law $\chi_\Phi(y) = \chi_{\Phi\setminus e}(y) - \chi_{\Phi/e}(y)$ for nonloops $e$ \cite{BG3}.  
The set-coloring function $\chiset_\D(k)$ is not a polynomial in $k$, of course, but also is not a polynomial function of $k! = |\fS_k|$ (see the small examples) and does not obey deletion-contraction for nonloops, not even with coefficients depending on $k$, as I found by computations with very small graphs.  A calculation with $\D = K_3$ convinced me the set-coloring function cannot obey deletion-contraction even if restricted to edges that are neither loops nor isthmi; but a second example would have to be computed to get a firm conclusion.  
However, going to gain graphs changes the picture: then there is a simple deletion-contraction law.  This indicates that the natural domain for studying set coloring and coloring from a partition is that of gain graphs.  I will develop this thought elsewhere.

%%%%%%%%%%%%%%%%%%%%%%%%

 \end{document}